\documentclass{article}

\usepackage{amssymb}

\usepackage[english,francais]{babel}

\newtheorem{theorem}{Theorem}[section]
\newtheorem{lemma}[theorem]{Lemma}
\newtheorem{e-proposition}[theorem]{Proposition}
\newtheorem{corollary}[theorem]{Corollary}
\newtheorem{e-definition}[theorem]{Definition\rm}

\def\og{\leavevmode\raise.3ex\hbox{$\scriptscriptstyle\langle\!\langle$~}}
\def\fg{\leavevmode\raise.3ex\hbox{~$\!\scriptscriptstyle\,\rangle\!\rangle$}}

\begin{document}
\title{{\bf Relative modular classes of Lie algebroids}}

\author{Yvette Kosmann-Schwarzbach
\thanks{UMR 7640 du CNRS}\\
Centre de Math\'ematiques Laurent 
Schwartz\\
\'Ecole Polytechnique\\
91128 Palaiseau, France\\
{\tt yks@math.polytechnique.fr}
\\
\\
Alan Weinstein
\thanks{Research  partially supported by NSF Grant
DMS-0204100}\\
Department of Mathematics\\
University of California\\Berkeley, CA
94720 USA\\
{\tt alanw@math.berkeley.edu}  
}
\date{}
\maketitle

\begin{abstract}
\noindent{\bf Abstract.}
We study the relative modular classes of Lie algebroids, and 
we determine their relationship with the 
modular classes of Lie algebroids with a twisted Poisson
structure.

\bigskip
\noindent{\bf R\'esum\'e. }
Nous \'etudions les classes modulaires relatives 
des alg\'ebro\"{\i}des de Lie
et nous d\'eterminons leur relation avec les
classes modulaires des alg\'ebro\"{\i}des de Lie avec structure de
Poisson tordue.
\end{abstract}

\section*{Version fran\c{c}aise abr\'eg\'ee}

Dans l'\'etude des structures de
Poisson sur les vari\'et\'es, Koszul a introduit 
un champ de vecteurs dont la diff\'erentielle
dans la cohomologie de Lichnerowicz-Poisson est nulle \cite{K}.
Dans \cite{W}, il est montr\'e que le flot d'un tel champ de vecteurs
est une limite classique du 
groupe d'automorphismes modulaires d'une alg\`ebre de von Neumann,
et la classe de cohomologie ainsi d\'efinie est 
nomm\'ee la {\it classe modulaire de la vari\'et\'e de Poisson}.
Il y est aussi indiqu\'e que cette notion se g\'en\'eralise 
lorsque le fibr\'e tangent \`a la vari\'et\'e est remplac\'e par
une alg\'ebro\"{\i}de de Lie arbitraire.
La notion de {\it classe caract\'eristique d'une
alg\'ebro\"{\i}de de Lie munie d'une repr\'esentation dans un fibr\'e
en droites} est d\'efinie dans \cite{ELW}. La {\it classe modulaire d'une
alg\'ebro\"{\i}de de Lie}, $A$, de base $M$ et d'ancre $\rho_A$,
est ensuite d\'efinie comme la classe caract\'eristique de $A$ pour la
repre\'sentation d\'efinie par (\ref{reprD}) avec $E=A$, $F =TM$,
$\varphi = \rho_A$.
La classe modulaire d'une vari\'et\'e de Poisson
en est un cas particulier 
lorsque l'alg\'ebro\"{\i}de de Lie consid\'er\'ee est le fibr\'e
cotangent de la vari\'et\'e muni du crochet de Lie des sections
d\'efini par le bivecteur de Poisson. Plus pr\'ecis\'ement, la classe 
modulaire du fibr\'e cotangent est le double de la classe modulaire de
la vari\'et\'e de Poisson. Ce facteur~2 s'explique par le fait
que, dans le cas g\'en\'eral, on doit consid\'erer un
fibr\'e en droites qui, dans le cas particulier, est le carr\'e 
de $\wedge^{\rm top}T^*M$. 
Le facteur 2 se retrouve dans les g\'en\'eralisations 
d\'ecrites ci-dessous.

Le probl\`eme de la g\'en\'eralisation au cas o\`u la structure de Poisson 
d'une vari\'et\'e est tordue par une 3-forme ferm\'ee, situation
\'etudi\'ee dans \cite{SW}, 
a \'et\'e r\'esolu dans \cite{KL}, o\`u 
une d\'efinition de la {\it classe modulaire d'une 
alg\'ebro\"{\i}de de Lie munie d'une structure de Poisson tordue} est
propos\'ee, le cas des vari\'et\'es de Poisson tordues 
\'etant celui o\`u
l'alg\'ebro\"{\i}de de Lie est le fibr\'e tangent \`a une
vari\'et\'e. 
Cette d\'efinition g\'en\'erale \'etend la caract\'erisation  donn\'ee
dans \cite{yks} (voir aussi \cite{Xu})~: 
le produit int\'erieur par un champ modulaire
est la diff\'erence de deux g\'en\'erateurs de carr\'e
nul de l'alg\`ebre de Gerstenhaber de l'alg\'ebro\"{\i}de de Lie. 

Dans cette Note, nous \'etudions d'abord la notion de 
{\it classe modulaire relative},
d\'efinie par la donn\'ee d'un morphisme d'alg\'ebro\"{\i}des de Lie,
notion introduite dans \cite{GMM} sous le nom de {\it classe modulaire d'un
morphisme d'alg\'ebro\"{\i}des de Lie}. Nous montrons ensuite
que la classe modulaire d\'efinie dans \cite{KL}
est la classe modulaire relative associ\'ee au morphisme 
$\pi^{\sharp} : A^* \to A$, d\'efini
par la donn\'ee d'une section $\pi$ de $\wedge^2 A$. Ici $A$ est une 
alg\'ebro\"{\i}de de Lie munie du bivecteur $\pi$ 
et d'une section $\psi$ de $\wedge^3(A^*)$ telle que
$d_A\psi =0$,  et satisfaisant la relation (\ref{eqtwisted});
l'espace des sections du 
fibr\'e vectoriel dual $A^*$ est muni du crochet de Lie associ\'e \`a 
la donn\'ee de $\pi$ et $\psi$.
Ce crochet, not\'e $[~,~]_{\pi,\psi}$, g\'en\'eralisant le
crochet de Lie des 1-formes sur une vari\'et\'e de 
Poisson, a \'et\'e d\'efini dans \cite{SW} lorsque $A$ est le fibr\'e
tangent \`a une vari\'et\'e de Poisson tordue, et \'etendu au
cas g\'en\'eral dans \cite{R}. La classe modulaire relative est une classe de
cohomologie de l'alg\'ebro\"{\i}de de Lie $A^*$, dont la
diff\'erentielle est not\'ee $d_{\pi,\psi}$.
(Pour une alg\'ebro\"{\i}de de Lie $E$, 
munie du crochet des sections not\'e en g\'en\'eral  $[~,~]_E$, 
on d\'esigne en g\'en\'eral par   
$d_E : \Gamma(\wedge^{\bullet}(E^*)) \to \Gamma(\wedge^{\bullet +
  1}(E^*))$ la diff\'erentielle de la cohomologie 
d'alg\'ebro\"{\i}de de Lie de $E$.) 

Nous rappelons d'abord (paragraphe \ref{modE})
la d\'efinition de la classe modulaire d'une
alg\'ebro\"{\i}de de Lie $E$. L'\'equation (\ref{defmodcl}) d\'efinit une  
section $\xi_E$ du fibr\'e vectoriel dual $E^*$ qui est ferm\'ee pour
la diff\'erentielle $d_E$ et qui d\'efinit par cons\'equent une classe
de cohomologie de degr\'e $1$ de $E$, not\'ee ${\rm {Mod}} \, E$.
Nous donnons en exemples $E=TM$ (la classe est nulle), $E=T^*M$
pour une vari\'et\'e de Poisson $(M,\pi)$, et 
le cas des alg\`ebres de Lie consid\'er\'ees comme 
alg\'ebro\"{\i}des de Lie de base un point (la classe est le
caract\`ere modulaire infinit\'esimal). La th\'eorie des classes
modulaires des alg\`ebres de Lie-Rinehart est d\'evelopp\'ee dans \cite{H2}.

Au paragraphe \ref{relmod}, nous d\'efinissons 
la classe modulaire relative 
d'un couple d'alg\'ebro\"{\i}des de Lie $(E,F)$
d\'efinie par un morphisme $\varphi$ de $E$ dans $F$.
Puis, nous d\'emontrons le
th\'eor\`eme \ref{representation}, 
r\'esultat g\'en\'eral qui sera utilis\'e dans 
l'application au cas des structures de Poisson tordues~:
la classe modulaire ${\rm{Mod}}^{\varphi} (E,F)$ est la classe
caract\'eristique, au sens de \cite{ELW}, de $E$ pour la
repr\'esentation 
$D^{\varphi}$, d\'efinie par l'\'equation (\ref{reprD}), 
de $E$ dans le fibr\'e vectoriel 
$L^{E,F} = \wedge^{\rm {top}} E \otimes \wedge^{\rm {top}} F^*$.
La d\'emonstration utilise le lemme \ref{lemma}
qui \'etend, du calcul tensoriel
sur les vari\'et\'es au cas d'une alg\'ebro\"{\i}de de Lie quelconque
$E$, la propri\'et\'e de
commutation de la d\'erivation de Lie des sections de
$\wedge^{\bullet}(E \oplus E^*)$ avec les contractions. 
 
Au paragraphe \ref{twistedpoisson}, on consid\`ere une
alg\'ebro\"{\i}de de Lie $A$, 
munie d'une structure de Poisson tordue.
Dans \cite{KL}, on a montr\'e que dans le cas particulier o\`u 
$A=TM$, la classe modulaire du fibr\'e cotangent de $A$ est \'egale
au double de la classe modulaire de $(TM, \pi, \psi)$ et l'on a
remarqu\'e que ce r\'esultat ne s'\'etend pas au cas g\'en\'eral.
Au th\'eor\`eme 4.1
nous \'etablissons la relation 
entre la classe relative du couple d'alg\'ebro\"{\i}des de Lie
$(A^*,A)$ d\'efinie par le morphisme $\pi^{\sharp}$ 
et la classe modulaire de $(A,\pi,\psi)$~: la premi\`ere est le double
de la seconde. 
Dans le cas particulier o\`u $A=TM$, la classe relative se r\'eduit
\`a la classe modulaire de $T^*M$ et l'on retrouve (corollaire
\ref{twiceTM}) le r\'esultat de \cite{KL}.

Le cas des alg\`ebres de Lie est abord\'e au 
paragraphe \ref{algebras}. Lorsque $\mathfrak h$ 
est une sous-alg\`ebre de Lie d'une alg\`ebre de Lie 
$\mathfrak g$, la classe modulaire de 
$(\mathfrak h, \mathfrak g)$ est l'obstruction \`a
l'existence d'une mesure $G$-invariante sur l'epsace homog\`ene
$G/H$, $H$ et $G$
\'etant des groupes de Lie connexes d'alg\`ebres de Lie $\mathfrak h$
et $\mathfrak g$. 

\selectlanguage{english}

\section{Introduction}

A modular class for a Lie algebroid with a twisted Poisson structure
$(A,\pi,\psi)$
was introduced in \cite{KL}.  In this note, we determine the relation
between the modular class of the Lie algebroid $A^*$
and this newly-defined
class. 
The vector bundle $A^*$ becomes a Lie algebroid when equipped with the
bracket $[~,~]_{\pi,\psi}$ on sections (defined in \cite{SW} and
\cite{R}) and the anchor $\rho_A \circ \pi^{\sharp}$, where $\rho_A$
is the anchor of $A$ and $\pi^{\sharp}:A^*\to A$ is the bundle map
associated to $\pi \in \Gamma(\wedge^2 A)$. 
To this end, we investigate
{\it relative modular classes}, which were first considered in
\cite{GMM}, where they were called {\it
modular classes of Lie algebroid morphisms}. 
We show that the relative modular class of $(A^*,A)$ defined by the
morphism $\pi^{\sharp}$ is equal to twice the class defined in
\cite{KL}
(Theorem 4.1).
Since the relative modular class of a pair of Lie algebroids $(E,F)$ defined
by a morphism $\varphi : E \to F$ is the difference of the modular
class of $E$ and the image under the dual $\varphi ^*$ of $\varphi$
of the modular class of $F$, it reduces to the modular class of $E$
when $F=TM$. So Theorem 4.1 extends the well-known fact 
\cite{ELW} \cite{W} that, for a Poisson manifold $(M,\pi)$, 
the modular class of the Lie algebroid $T^*M$ is twice the modular
class of the Poisson manifold. 
To prove Theorem 4.1,
we apply Theorem 3.3,
which states that the relative modular
class of the pair $(E,F)$ defined by $\varphi$ is a characteristic class,
in the sense of \cite{ELW}, 
of $E$ equipped with the representation defined by (\ref{reprD}).

\section{The modular class of a Lie algebroid}\label{modE}

The {\it modular class} ${\rm Mod}\,{E}$ 
of a Lie algebroid
$(E,\rho_{E})$, with base $M$, was defined by Evens, Lu and
Weinstein in \cite{ELW}. 
It is the class of a section $\xi_E$ of $E^{*}$ which
satisfies, for all $x \in \Gamma E$,
\begin{equation}\label{defmodcl}
< \xi _{E} , x> \, \omega \otimes \lambda = [x,
\omega]_{E} \otimes \lambda + \omega \otimes
\mathcal{L}_{\rho_{E}x}  \lambda \ ,
\end{equation}
where $\omega$ (resp., $\lambda$) is a nowhere-vanishing section
of $\wedge^{\rm{top}}E$ (resp.,  $\wedge^{\rm{top}} T^{*}M$).
Since the definition of $\xi_E$ is insensitive to changes in the sign
of $\omega$ and $\lambda$, we may use densities rather than forms if
$E$ or $T^*M$ is non-orientable.  
The section $\xi_{E}$ is called a {\it modular section} of $E$. 
It is a 1-cocycle in the Lie algebroid cohomology
of $E$ and its class,  ${\rm Mod}\,E$, 
is independent of the choice of $\omega$ and
$\lambda$. 
We list various examples. 

(i) Let $M$ be a manifold and $E= TM$, with anchor the identity of
$TM$. Then ${\rm Mod}(TM) = 0$.

(ii) Let $(M,\pi)$ be a Poisson manifold and $E= T^{*}M$ with anchor
$\pi^{\sharp}$ and Lie bracket $[\ ,\ ]_{\pi}$. Then the class of
${\rm Mod}(T^{*} M)$ is twice the modular class of $(M,\pi)$, as defined in
\cite{W}, following the earlier definition of the modular vector
field, without a name, in \cite{K},
and its use in \cite{DH}, \cite{G} and \cite{LX}. 

(iii) If ${\mathfrak g}$ is a Lie algebra, both the modular class of 
${\mathfrak g}$ considered as a Lie algebroid with
base a point and the modular class of the
linear Poisson manifold $\mathfrak{g}^{*}$ are equal to 
the {\it infinitesimal modular character} 
 $\chi^{\mathfrak g}$ of ${\mathfrak g}$, which is  
the linear form 
$x \mapsto {\rm{Tr}}({\rm{ad}}^{\mathfrak g}_x)$ on $\mathfrak g$.

\section{Relative modular classes}\label{relmod}

Let $(E,\rho_{E})$ and $(F,\rho_{F})$ be Lie algebroids with base $M$.
We recall that a vector bundle morphism $\varphi:E\to F$
is a Lie algebroid morphism (over the identity on $M$) if and only if 
$\wedge^{\bullet} \varphi^{*}$ is a chain map from
$(\Gamma(\wedge^{\bullet}(F^*), d_F)$ to
$(\Gamma(\wedge^{\bullet}(E^*), d_E)$.

\begin{e-definition}\label{defrelative}For a Lie algebroid morphism
  $\varphi:E\to F$, we set
\begin{equation}\label{def}
{\rm{Mod}} ^{\varphi} (E,F) = 
{\rm{Mod}} \, E - \varphi^{*} {\rm{Mod}} \, F \ ,
\end{equation}
and we call ${\rm{Mod}}
^{\varphi} (E,F)$
the {\em relative modular class} of $(E,F)$ defined by $\varphi$.
\end{e-definition}

For Lie algebroids $A, B, C$ and morphisms 
$\varphi : A \to B$ 
and $\psi : B \to C$, we obtain
\begin{equation}
{\rm{Mod}}^{\psi\circ\varphi} (A,C) = {\rm{Mod}} ^{\varphi} (A,B)
+\varphi^{*} {\rm{Mod}}^{\psi} (B,C)  \ .
\end{equation}  

It is clear from the definition that
\begin{equation}\label{composition}
{\rm{Mod}}^{\rho_{E}} (E,TM) = {\rm{Mod}} \, E \ ,
\end{equation} 
so the ``absolute'' modular class is a special case of the
 relative  class. 
 
This statement was a theorem in \cite{GMM}, where the
relative class was defined in terms of divergence operators,
without the use of the absolute class.  
The proof of that theorem is related to the following argument, in which we
show directly that the relative modular class
${\rm Mod}^\varphi(E,F)$ is the characteristic class of $E$
associated with a representation, a notion defined in \cite{ELW}.

We shall first prove a lemma.
Let $E$ be a Lie algebroid and let $x$ be a section of $\Gamma E$.  
On the one hand, the Lie derivative with respect to $x$ of a form $\alpha
\in \Gamma (\wedge^{\bullet} E^{*})$ is 
$\mathcal{L}^{E}_{x}\alpha = [i_{x},d_{E}]\alpha$, and  
$\mathcal{L}^{E}_{x}$ is a  degree $0$ derivation 
of $\Gamma (\wedge^{\bullet} E^{*})$. On the other hand, the map $Q
 \mapsto [x,Q]_{E}$ from $\Gamma(\wedge^{\bullet} E)$ to
$ \Gamma (\wedge^{\bullet} E)$
is a degree $0$ derivation of $\Gamma (\wedge^{\bullet} E)$. 
If $f \in C^{\infty}(M)$ is considered as a 0-form
or a 0-vector, then $\mathcal{L}^{E}_{x}f = <d_Ef,x> = \rho_E (x)\cdot f = 
[x,f]_E$. Therefore these operations extend to a unique derivation of 
$\Gamma (\wedge^{\bullet} (E \oplus E^{*}))$,
which we also denote by $\mathcal{L}^{E}_{x}$.

\begin{lemma} \label{lemma}
For $x\in \Gamma E$ and $\alpha \in \Gamma
  (\wedge^{\bullet}E^{*})$,
the endomorphisms $\mathcal{L}^{E}_{x}$, $i_{\alpha}$ and 
$i_{\mathcal{L}^{E}_{x}\alpha}$ of  
$\Gamma (\wedge ^{\bullet} E)$ satisfy
\begin{equation}\label{commutation} 
[\mathcal{L}^{E}_{x} , i_{\alpha}] = i_{\mathcal{L}^{E}_{x}\alpha} \ .
\end{equation}
\end{lemma}

\noindent{\it Proof.}
If $\alpha$ is a 0-form, the relation follows from the Leibniz rule
for $[~,~]_E$.
Next, let $\alpha$ be a 1-form. Since the derivation $
[\mathcal{L}^{E}_{x} , i_{\alpha}] - i_{\mathcal{L}^{E}_{x}\alpha}$
vanishes on elements of degree $0$, it is enough to prove that it
vanishes on elements of degree $1$. Let  $y$ be a section of $E$, then
$[\mathcal{L}^{E}_{x} , i_{\alpha}](y) = \rho_E(x) \, \cdot
<\alpha, y> - <\alpha, [x,y]_E>$,
while 
$ 
i_{\mathcal{L}^{E}_{x}\alpha} y = <{\mathcal{L}^{E}_{x}\alpha}, y> = 
 \rho_E(y) \cdot
<\alpha, x> + (d_E\alpha) (x,y)$.
By the definition of $d_E$ the difference 
$[\mathcal{L}^{E}_{x} , i_{\alpha}] - i_{\mathcal{L}^{E}_{x}\alpha}$
vanishes on $y$.
To extend the formula to forms of arbitrary degree, we use
induction on the degree, taking into account  
 that $i_{\alpha \wedge \alpha_1}
= i_{\alpha} \circ i_{\alpha_1}$ for forms $\alpha$ and
$\alpha_1$ and 
that $\mathcal{L}^{E}_{x}$ is a  degree $0$ derivation of
$\Gamma(\wedge^{\bullet} E^*)$.
\hfill $\square$

\begin{theorem}\label{representation} 
Let $\varphi \, : \, E \to F$ be a Lie algebroid
  morphism. 
Set $L^{E,F} = \wedge^{\rm {top}} E \otimes \wedge^{\rm {top}} F^*$, and let
$\omega \otimes \nu \in \Gamma(L^{E,F})$.
For $x\in \Gamma E$, set 
\begin{equation}\label{reprD}
D^{\varphi}_{x} (\omega\otimes \nu) = \mathcal{L}^{E}_{x}\omega 
\otimes \nu
+ \omega \otimes \mathcal{L}^{F}_{\varphi x} \nu \ .
\end{equation}
\noindent(i)
The map  $x\mapsto D_{x}^{\varphi}$ is a representation of $E$ on
$L^{E,F} $. \\
\noindent(ii)
The relative modular class ${{\rm{Mod}}}^\varphi (E,F)$ 
is the characteristic class of the Lie algebroid $E$ with
representation $D^{\varphi}$.
\end{theorem}

\noindent{\it Proof.}
(i) We must prove that 
$$
{\rm {(a)}} \,
D^{\varphi}_{fx} = fD_{x}^{\varphi} \ , \quad 
{\rm {(b)}} \, D^{\varphi}
_{x} (f(\omega \otimes \nu))= f D^{\varphi}_{x} (\omega \otimes\nu)+ (\rho_{E}
(x) \cdot f) \, \omega \otimes \nu \ ,
$$
and
$$ 
{\rm {(c)}} \, D^{\varphi}_{[x,y]_E} = [D^{\varphi}_{x} ,D^{\varphi}_{y}] \ , 
$$ 
for all $f \in C^{\infty}(M)$, $x$ and $y \in \Gamma E$.
In fact, since $[f,~ .~]_E$ and $i_{d_E f}$ are derivations of 
$\Gamma(\wedge^{\bullet}E)$ of degree $-1$ 
which coincide on elements of degree 0 or 1, 
and are therefore equal,
$
[fx,\omega]_E- f[x,\omega]_E= [f,\omega]_E \wedge x =i_{d_Ef}\omega 
\wedge x = (\rho_E(x) \cdot f) \omega$.
In addition,
$
d_{F}i_{f\varphi x} \nu - f d_Fi_{\varphi x} \nu =
 d_{F} f \wedge i_{\varphi x} \nu = 
< \varphi x, d_{F} f > \nu = (\rho_{E} (x) \cdot f) \nu \ ,
$
since
$\rho_{F} \circ \varphi =\rho_{E}$. 
Now (a) follows from these equalities.
Relation (b)
is clear, and 
(c) follows from the morphism property of
$\varphi$.

(ii) 
We shall use relation (\ref{defmodcl})  for $\xi_{E}$ and $\xi_{F}$, and 
Lemma \ref{lemma}. Let $\mu \in \Gamma (\wedge^{\rm{top}}
T^{*}M)$ be a volume form on $M$, $\omega \in \Gamma
(\wedge^{\rm{top}}E)$, and $\bar \omega \in \Gamma
(\wedge^{\rm{top}}F)$, both nowhere-vanishing.
Thus $\xi_{E}$ satisfies, for all $x\in \Gamma E$,
\begin{equation}\label{xiE}
< \xi_{E} , x > \, \omega\otimes \mu = \mathcal{L}^{E}_{x} \omega \otimes \mu
+\omega \otimes \mathcal{L}_{\rho_{E}x}\mu \ ,
\end{equation}
while $\xi_{F}$
satisfies, for all $x\in \Gamma E$,
\begin{equation}\label{xiF}
< \varphi^{*} \xi_{F} , x > \, \bar\omega \otimes \mu =
 \mathcal{L}^{F}_{\varphi x}\bar \omega \otimes \mu + \bar\omega \otimes
\mathcal{L}_{\rho_{F} ( \varphi x)} \mu \ .
\end{equation}
Let $\nu$ be the section of $\wedge^{\rm{top}}F^{*}$ such that
$< \bar\omega ,\nu > = 1$.
By Lemma \ref{lemma}, 
$ < {\mathcal L}_{\varphi x}^F \bar\omega , \nu > 
+ < \bar \omega , \mathcal{L}^{F}_{\varphi x} \nu > = 0$. 
Setting $\eta = \xi_{E} -\varphi^{*}\xi_{F}$, we
obtain, for all $x \in \Gamma E$,
\begin{equation}\label{repr}
< \eta , x > \, \omega \otimes \nu =
{\mathcal L}^E_x \omega \otimes \nu + \omega \otimes
\mathcal{L}^{F}_{\varphi x} \nu \ .
\end{equation}
Therefore, the section $\eta$ of $E^*$ is a representative of
the class ${{\rm{Mod}}} ^{\varphi} (E,F)$. 
\hfill $\square$

\section{Lie algebroids with  twisted Poisson
  structures}\label{twistedpoisson}

A Lie algebroid $A$  equipped with 
a bivector $\pi \in \Gamma (\wedge^{2}A)$ and a $d_A$-closed 3-form
$\psi\in \Gamma (\wedge^{3} A^{*})$ satisfying 
\begin{equation}\label{eqtwisted}
\frac{1}{2} [\pi ,\pi]_{A} = (\wedge^{3}\pi^{\sharp})\psi \ ,
\end{equation}
is called a {\it Lie algebroid with a twisted 
Poisson structure} \cite{KL}. 
The vector bundle $A^{*}$ is then a Lie algebroid with anchor
$\pi^{\sharp}$ and Lie bracket of sections $[\  ,\ ]_{A^{*}} = [\ ,\
]_{\pi ,\psi}$, as in \cite{SW}.
Such structures on general Lie algebroids were first studied in \cite{R}. 
The untwisted case, $\psi = 0$, is that of Poisson structures on Lie
algebroids; the pair $(A,A^{*})$ is then a Lie
bialgebroid, which is called {\it triangular}.

If $A=TM$, then 
$(M, \pi ,\psi)$ is called 
a {\it twisted Poisson manifold}, for which see \cite{SW}.
The case when $A=TM$ and $\psi = 0$ is that of Poisson manifolds.

When $(A,\pi,\psi)$ is a Lie algebroid 
with a twisted Poisson structure, 
we may compare 
two cohomology classes of $A^{*}$: the modular class ${\rm{Mod}}(A^*)$
and the modular class defined in \cite{KL}
which we shall denote by $\theta_{KL}(A,\pi,\psi)$.

\begin{theorem} \label{twisted}
The cohomology classes  
${\rm {Mod}}(A^*)$ and $\theta_{KL}(A,\pi,\psi)$ 
are related by  
\begin{equation}\label{twice}
2\theta_{KL}(A,\pi,\psi) = {{\rm{Mod}}}^{\pi^{\sharp}}(A^{*}, A) =
  {\rm{Mod}}(A^*) - (\pi^{\sharp})^* ({\rm{Mod}}A) \ .
\end{equation}  
\end{theorem}
\noindent{\it Proof.} Let $\lambda$ be a nowhere-vanishing 
section of $\wedge^{\rm{top}}
  A^{*}$.
By Theorem \ref{representation},
 a representative $W\in \Gamma A$ of the relative modular
class ${{\rm{Mod}}}^{\pi^{\sharp}}(A^{*},A)$ satisfies, for all $\alpha \in
\Gamma A^{*}$,
\begin{equation}\label{eqW}
<\alpha, W > \, \lambda \otimes \lambda = [\alpha ,
\lambda]_{\pi ,\psi} \otimes \lambda + \lambda \otimes
\mathcal{L}^{A}_{\pi  ^{\sharp} \alpha} \lambda\ .
\end{equation}
Adopting the notations of \cite{KL}, we set $\partial_{\pi}= [d_A,
i_{\pi}]$ and $Y_{\pi
  ,\psi} = \pi^{\sharp} i_{\pi} \psi$,
and we consider the generator 
$\partial = \partial_{\pi} + \underline{\partial}_{\pi,\psi} +
i_{Y_{\pi,\psi}}$
of $[\ ,\ ]_{\pi ,\psi}$,
which is of square $0$. (Here ${\underline
  {\partial}}_{\pi, \psi}$ is an operator on forms, which vanishes on functions
and $1$-forms.)
Since $\lambda$ is of top degree, $\partial$ satisfies the relation
$
[\alpha , \lambda]_{\pi ,\psi} = (\partial \alpha)\lambda - \alpha
\wedge \partial \lambda$ for all $\alpha \in \Gamma (A^*)$.

The class $\theta_{KL}(A,\pi,\phi)$ is represented by
the modular section $Z_{\pi,\psi,\lambda} = X_{\pi,\lambda} +
Y_{\pi, \psi}$,
where $X_{\pi,\lambda}$ satisfies
$
\mathcal{L}^{A}_{\pi^{\sharp}\alpha} \lambda =  <\alpha , X_{\pi,\lambda}
> \lambda - (\partial_{\pi} \alpha) \lambda$, and $Z_{\pi,\psi,\lambda}$
satisfies 
$\partial\lambda = -i_{Z_{\pi,\psi,\lambda} } \lambda$ (see (6.8) in
\cite{KL}).
Therefore 
$<\alpha, W > \, \lambda \otimes \lambda = (\partial \alpha) \lambda
\otimes \lambda + (\alpha \wedge i_{Z_{\pi,\psi,\lambda} } \lambda ) 
\otimes \lambda
+ < \alpha, X_{\pi, \lambda}> \lambda \otimes \lambda - (\partial
_{\pi} \alpha) \lambda \otimes \lambda$. Since
$ \partial \alpha - 
\partial_{\pi} \alpha = \\
<\alpha , Y_{\pi ,\psi} >$, we obtain
$< \alpha , W > \lambda \otimes \lambda 
= < \alpha , Y_{\pi,\psi} > \lambda \otimes \lambda 
+ < \alpha , Z_{\pi,\psi,\lambda}> \lambda \otimes \lambda 
+ < \alpha , X_{\pi,\lambda}> \lambda \otimes \lambda$,
whence
$ W = 2 Z_{\pi,\psi,\lambda}$.
\hfill $\square$

It was proved in \cite{KL} that 
$\theta_{KL}(A,\pi, \psi)$ is the characteristic
class of the Lie algebroid $A^*$ with representation $D^{\partial}$ in
$\wedge^{\rm{top}}(A^*)$, where 
$D^{\partial}_{\alpha} \lambda = - \alpha \wedge \partial
\lambda$. It is easy to show that 
$D^{\pi^{\sharp}}(\lambda \otimes \lambda) = D^{\partial} \lambda
\otimes \lambda + \lambda \otimes
D^{\partial} \lambda$, i. e., $D^{\pi^{\sharp}}$ is the ``square'' of
$D^{\partial}$. 
This remark leads to an alternate proof of Theorem \ref{twisted}.

Define the {\it modular class of a twisted Poisson manifold} $(M,\pi,\psi)$
to be $\theta_{KL}(TM,\pi,\psi)$.
Applying Theorem \ref{twisted} to $A= TM$ we recover the following result
of \cite{KL}.
\begin{corollary}\label{twiceTM} 
Let $(M,\pi,\psi)$ be a twisted Poisson manifold. Then
$2 \theta_{KL}(TM,\pi,\psi) = {\rm {Mod}}(T^*M)$. 
\end{corollary}

When both $A=TM$ and $\psi =0$, the class $\theta_{KL}(TM,
\pi,\psi)$ reduces to the modular class of the Poisson manifold
$(M,\pi)$. Thus, 
this corollary generalizes to the twisted case the property of the 
modular class of a Poisson manifold that we recalled in Section \ref{modE}.

\section{The case of Lie algebras}\label{algebras}

If ${\mathfrak a}$ and  ${\mathfrak b}$ are Lie algebras, they can be
considered as Lie algebroids with base a point. Thus, if $\varphi : 
{\mathfrak a} \to {\mathfrak b}$ is a homomorphism of Lie algebras,
${\rm {Mod}}^{\varphi}({\mathfrak a},{\mathfrak b})$ is a class in the Lie
algebroid cohomology of ${\mathfrak a}$ of degree~$1$, i. e., a
1-cocycle in ${\mathfrak a}^*$.

When ${\mathfrak h}$ is a Lie subalgebra of ${\mathfrak g}$ and
$\varphi$ is the canonical injection $\iota$, then 
${\rm {Mod}}^{\iota}({\mathfrak h},{\mathfrak g})= 
\chi^{\mathfrak h} - \iota^*(\chi^{\mathfrak g})$.
When $H$ is a connected closed Lie subgroup with Lie algebra 
${\mathfrak h}$ of a connected Lie
group $G$ with Lie algebra ${\mathfrak g}$, 
the vanishing of the relative modular class 
${\rm {Mod}}^{\iota}({\mathfrak h},{\mathfrak g})$ is  necessary and
sufficient for the existence of a $G$-invariant measure on the
homogeneous space $G/H$. 
This follows from the fact, proved in \cite{Weil},
that there exists a $G$-invariant measure on $G/H$ if and only if 
$\Delta^H - \iota^* \Delta^G = 0$, where $\Delta^G$ (resp., $\Delta^H$)
is the modular function of $G$ (resp., $H$), and the 
relation ${\rm {Det}}( {\rm {Ad}}^G_{\exp tx}) = 
\exp (t \chi^{\mathfrak g} (x))$.

\end{document}